\documentclass[hyper,11pt]{JHEP3}

\usepackage{xy,amssymb,amsmath,diagrams,times,theorem,longtable}

\xyoption{poly}
\xyoption{arc}
\xyoption{knot}
\xyoption{curve}
\xyoption{arrow}
\xyoption{all}

\newcommand{\C}{\mathbb{C}}

\newcommand{\N}{\mathbb{N}}

\newcommand{\Z}{\mathbb{Z}}

\newcommand{\Oscr}{{\mathcal O}}

\newcommand{\Af}{{\mathbb A}}

\newcommand{\vtx}[1]{*+[o][F-]{\scriptscriptstyle #1}}

\newcommand{\wis}[1]{{\text{\em \usefont{OT1}{cmtt}{m}{n} #1}}}

\newtheorem{theorem}{Theorem}

\newtheorem{lemma}{Lemma}

{\theorembodyfont{\rmfamily} 
\newtheorem{example}{Example}

\newtheorem{definition}{Definition} }

\preprint{UIA preprint 2003-04}

\title{One quiver to rule them all}

\author{Lieven Le Bruyn \\
Departement Wiskunde en Informatica, Universiteit Antwerpen  \\
B-2020 Antwerp (Belgium) \\
E-mail :  \email{lieven.lebruyn@ua.ac.be}}

\abstract{In math.AG/9904171 it was shown that the \'etale local structure of finite dimensional representations for a formally smooth algebra is determined by (varying) local quiver settings. In this note we prove that there is one quiver setting $(Q_A,\alpha_A)$ depending only on the formally smooth algebra $A$ which contains enough information to reconstruct all these local quiver settings. Conjecturally, the formally smooth algebra $A$ is locally isomorphic in a (yet to be developed) non-commutative \'etale topology to an algebra $B$ Morita equivalent (determined by the dimension vector $\alpha_A$) to the path algebra $\C Q_A$.}

\begin{document}

\begin{center}
{\it In memory of International Law $\quad \dagger \quad$ 20-03-2003}
\end{center}

\vskip 5mm

Following \cite{Kontsevich} and \cite{CuntzQuillen} one defines a non-commutative smooth affine variety to correspond to a {\em formally smooth algebra} $A$. Such an algebra $A$ has the lifting property with respect to nilpotent ideals in $\wis{alg}$, the category of all associative $\C$-algebras with unit. That is, for every $B \in \wis{alg}$, every nilpotent ideal $I \triangleleft B$ and every $\C$-algebra morphism $\phi : A \rTo B/I$, there exists a lifted algebra morphism $\tilde{\phi}$
making the diagram below commutative
\[
\begin{diagram} 
A & \rDotsto^{\tilde{\phi}} & B \\
& \rdTo_{\phi} & \dOnto \\
& & \frac{B}{I}
\end{diagram}
\]
This notion generalizes Grothendieck's characterization of commutative regular algebras (replacing $\wis{alg}$ by $\wis{commalg}$, the category of all commutative $\C$-algebras) but as the lifting property in $\wis{alg}$ is stronger not all commutative regular algebras will be formally smooth. In fact, by \cite{CuntzQuillen} any commutative formally smooth affine algebra is the coordinate ring of a disjoint union of points and smooth affine curves.

Typical non-commutative examples of formally smooth algebras are path algebras $\C Q$ of finite quivers $Q$ (see \cite{LBnagatn}), in particular free associative algebras $\C \langle x_1,\hdots,x_n \rangle$. Following \cite{Kontsevich} (or \cite{LBnagatn}) one assigns to a formally smooth affine algebra $A$ the family of finite dimensional representation schemes $\{ \wis{rep}_n~A : n \in \N \}$ each element of which is a smooth affine (commutative) variety (possibly containing several connected components). There is a natural base-change action by $GL_n$ on $\wis{rep}_n~A$ with quotient variety $\wis{iss}_n~A$ parametrizing isomorphism classes of semi-simple $n$-dimensional representations. The \'etale local structure of the quotient varieties $\wis{iss}_n~A$ is described in terms of local quiver-settings which we will recall briefly.

\section{Local quiver settings}

A point $\xi$ in $\wis{iss}_n~A$ corresponds to the isomorphism class of a semi-simple $n$-dimensional representation of $A$
\[
M_{\xi} = S_1^{\oplus e_1} \oplus \hdots \oplus S_k^{\oplus e_k} \]
where the $S_i$ are the distinct simple components (say of dimension $d_i$) occurring in $A$ with multiplicity $e_i$. Construct a quiver $Q_{\xi}$ on $k$-vertices $\{ v_1,\hdots,v_k \}$ (corresponding to the distinct simple components of $M_{\xi}$) with the property that the number of directed arrows from $v_i$ to $v_j$ is equal to the dimension of the extension group $Ext^1_A(S_i,S_j)$. Consider the dimension vector $\alpha_{\xi} = (e_1,\hdots,e_k)$ (corresponding to the multiplicities of the simple components in $M_{\xi}$) and recall that $\wis{rep}_{\alpha_{xi}}~Q_{\xi}$ is the affine space of all $\alpha_{\xi}$-dimensional representations of the quiver $Q_{\xi}$. On this space there is a base-change action by the group
$GL(\alpha_{\xi}) = GL_{e_1} \times \hdots \times GL_{e_k}$ and the corresponding affine quotient variety $\wis{iss}_{\alpha_{\xi}}~Q_{\xi}$ parametrizes semi-simple $\alpha_{\xi}$-dimensional representations of $Q_{\xi}$, see \cite{LBProcesi}. As $n = \sum_i e_id_i$ there is a natural embedding of $GL(\alpha_{\xi})$ in $GL_n$. With these notations, the \'etale local structure of the quotient variety $\wis{iss}_n~A$ near $\xi$ is described by the following result of \cite{LBnagatn}.

\begin{theorem} If $A$ is a formally smooth affine algebra, there is a $GL_n$-equivariant \'etale isomorphism between
\[
\wis{rep}_n~A \qquad \text{and} \qquad GL_n \times^{GL(\alpha_{\xi})} \wis{rep}_{\alpha_{\xi}}~Q_{\xi} \]
in a neighborhood of the orbit of $M_{\xi}$ (resp. the orbit of the zero representation). As a consequence, there is an \'etale isomorphism between
\[
\wis{iss}_n~A \qquad \text{and} \qquad \wis{iss}_{\alpha_{\xi}}~Q_{\xi} \]
in a neighborhood of $\xi$ (resp. a neighborhood of the point $\overline{0}$ corresponding to the zero representation).
\end{theorem}

The main result of this note asserts that there is one quiver setting $(Q_A,\alpha_A)$ depending only on the formally smooth algebra $A$ that contains enough information to reconstruct all these local quiver settings $(Q_{\xi},\alpha_{\xi})$ for $\xi \in \wis{iss}_n~A$ for any $n \in \N$.

\section{The quiver setting $(Q_A,\alpha_A)$}

If $A$ is a formally smooth affine algebra we can decompose the affine smooth variety
\[
\wis{rep}_n~A = \bigsqcup_{| \alpha | = n}~\wis{rep}_{\alpha}~A \]
into its connected components $\wis{rep}_{\alpha}~A$ and we will call the label $\alpha$ a dimension vector of $A$ of total dimension $n$. The set $\wis{comp}~A$ of all dimension vectors of $A$ can be equipped with an Abelian semigroup structure by defining $\alpha+\beta = \gamma$ whenever $\wis{rep}_{\gamma}~A$ is the connected component of $\wis{rep}_{m+n}~A$ containing the image of $\wis{rep}_{\alpha}~A \times \wis{rep}_{\beta}~A$ under the direct sum map
\[
\wis{rep}_n~A \times \wis{rep}_m~A \rTo^{\oplus} \wis{rep}_{m+n}~A \]
for $\alpha$ a dimension vector of total dimension $n$ and $\beta$ of total dimension $m$. The semigroup $\wis{comp}~A$ is augmented as there is a natural map $\wis{comp}~A \rTo \N$ sending $\alpha$ to its total dimension $| \alpha |$. In fact, for a formally smooth algebra $A$ this definition of $\wis{comp}~A$ coincides with the component semigroup introduced and studied by K. Morrison in \cite{Morrison}.

Clearly, $\wis{comp}~\C Q \simeq \N^k$ if $Q$ has $k$ vertices as the semigroup generators correspond to the vertex-simples. If $X$ is a smooth affine curve, then $\wis{comp}~\C[X] \simeq \N$ as all $\wis{rep}_n~\C[X]$ are irreducible smooth affine varieties. In general though the structure of $\wis{comp}~A$ can be quite complicated.

\begin{example} For every sub-semigroup $S \subset \N$ there exists an affine formally smooth algebra $A$ such that
\[
\wis{comp}~A \simeq S \]
as augmented semigroups.
\end{example}

\Proof
A special affine algebra of inversion depth $b$ $A$ has a presentation
\[
A = \frac{\C \langle x_1,\hdots,x_a,y_1,\hdots,y_b \rangle}{(1-y_ip_i(x_1,\hdots,x_a,y_1,\hdots,y_{i-1})~,~1 \leq i \leq b )}
\]
where the $p_i$ are polynomials in the noncommuting variables $\{ x_1,\hdots,x_a,y_1,\hdots,y_{i-1} \}$. Therefore, $A$ is a universal localization (see e.g. \cite{Schofield} for definition and properties) of the free associative algebra $\C \langle x_1,\hdots,x_a \rangle$ and as such $A$ is formally smooth. Because $\wis{rep}_n~\C \langle x_1,\hdots,x_a \rangle = M_n(\C)^{\oplus a}$ is irreducible, $\wis{comp}~A \subset \N$
and consists of those $n \in \N$ such that none of the $p_i$ is a rational identity for $n \times n$ matrices (see \cite{LBratid}).

In the special case when $gcd(S) = 1$ it was proved in \cite[Coroll. 1]{LBratid} that there is a special affine algebra $A$ such that $\wis{comp}~A = S$. In general, if $gcd(S) = n$ let $S' = S/n$ and $A'$ a special affine algebra such that $\wis{comp}~A' = S'$. Then, the free algebra product
$A = M_n(\C) \ast A'$ is again a formally smooth algebra and satisfies $\wis{comp}~A = S$.

\par \vskip 4mm
In all these examples, the component semigroup is finitely generated but we do not know whether this is always the case for an affine formally smooth algebra. From now on, we will impose :

\par \vskip 4mm
\noindent
{\bf Assumption 1 :} $A$ is an affine formally smooth algebra such that $\wis{comp}~A$ has a finite number of semigroup generators.

\par \vskip 4mm
There exists a ringtheoretic characterization of these semigroup generators :

\begin{lemma} The following are equivalent
\begin{enumerate}
\item{$\alpha$ is a semigroup generator of $\wis{comp}~A$ of total dimension $n$.}
\item{The quotient map $\wis{rep}_{\alpha}~A \rOnto \wis{iss}_{\alpha}~A$ is a principal $PGL_n$-fibration in the \'etale topology.}
\item{The ring $\int_{\alpha}~A$ of all $GL_n$-equivariant maps from $\wis{rep}_{\alpha}~A$ to $M_n(\C)$ is an Azumaya algebra with center $\C[\wis{iss}_{\alpha}~A]$.}
\end{enumerate}
\end{lemma}

\Proof
If $M \in \wis{rep}_{\alpha}~A$ is not a simple representation, then the orbit closure $\Oscr(M)$ contains a semi-simple representation corresponding to the Jordan-H\"older decomposition of $M$ into distinct simple components
\[
N = S_1^{\oplus e_1} \oplus \hdots \oplus S_l^{\oplus e_l} \]
where $S_i$ is simple of dimension vector $\beta_i$. But then, $\alpha = \sum_i e_i \beta_i$. As the stabilizer subgroup of a simple representation is $\C^*$ (by Schur's lemma) this proves that $1$ and $2$ are equivalent. The equivalence of $2$ and $3$ follows from \cite{Procesi} and the fact that both are classified by the \'etale cohomology group $H^1_{et}(\wis{iss}_{\alpha}~A,PGL_n)$.

\par \vskip 4mm
For $\alpha,\beta \in \wis{comp}~A$ and $N \in \wis{rep}_{\alpha}~A$. $M \in \wis{rep}_{\beta}~A$ we know from \cite[Lemma 4.3]{CBS} that the dimension of the extension space
\[
dim_{\C}~Ext^1_A(N,M) \]
is an upper semi-continuous function on $\wis{rep}_{\alpha}~A \times \wis{rep}_{\beta}~A$. In particular, there is a Zariski open subset $E(\alpha,\beta)$ of this product where this dimension attains its minimal value which we will denote by $ext(\alpha,\beta)$.

\begin{definition} Let $A$ be an affine formally smooth algebra having a finite set of semigroup generators $\{ \beta_1,\hdots,\beta_k \}$ of $\wis{comp}~A$. Let $Q_A$ be the quiver on $k$ vertices $\{ v_1,\hdots,v_k \}$ such that the number of directed arrows from $v_i$ to $v_j$ is equal to $ext(\beta_i,\beta_j)$ and the number of loops in $v_i$ is equal to $dim~\wis{iss}_{\beta_i}~A$. Let $\alpha_A$ be the dimension vector $(n_1,\hdots,n_k)$ where $n_i = | \beta_i |$. We say that $A$ is a formally smooth algebra of type $(Q_A,\alpha_A)$.
\end{definition}

\section{Some examples}

Let $A$ be an affine formally smooth algebra of type $(Q_A,\alpha_A)$ with $\alpha_A = (n_1,\hdots,n_k)$ and construct the algebra
\[
B = \begin{bmatrix}
B_{11} & \hdots & B_{1k} \\
\vdots & & \vdots \\
B_{k1} & \hdots & B_{kk}
\end{bmatrix}
\]
where $B_{ij}$ is the $n_i \times n_j$ block matrix with all its components equal to the subspace of
$\C Q_A$ spanned by all oriented paths in $Q_A$ starting at $v_i$ and ending at $v_j$. That is, $B$ is an affine formally smooth algebra which is Morita equivalent to the path algebra $\C Q_A$ with the Morita equivalence determined by the dimension vector $\alpha_A$. The hope is that in a (yet to be developed) non-commutative \'etale topology the algebras $A$ and $B$ are locally isomorphic. The examples below may add some weight to this conjecture.

\begin{example} Let $X$ be an affine smooth curve with coordinate ring $A = \C[X]$, then $\wis{comp}~A = \N$ with generating component $\wis{rep}_1~A = \wis{iss}_1~A = X$. Therefore,
$A$ is of type
\par \vskip 3mm
\[
Q_A = \xymatrix{\vtx{}\ar@(lu,ru)} \qquad \text{and} \qquad \alpha_A = (1) \]
The associated algebra $B$ is in this case isomorphic to the polynomial ring $\C[x]$ which is indeed locally isomorphic (in the \'etale topology) to $A$.
\end{example}

\begin{example} If $A$ is the path algebra $\C Q$ then $\wis{comp}~A = \N^k$ with semigroup generators the vertex dimension vectors $\delta_i = (0,\hdots,0,1,0,\hdots,0)$. Clearly,
\[
\wis{rep}_{\delta_i}~A = \wis{iss}_{\delta_i}~A = \Af^{l_i} \]
where $l_i$ is the number of loops in the $i$-th vertex of $Q$. If $\delta_i \not= \delta_j$ then we obtain from the Euler-form formula that $dim_{\C}~Ext^1_A(S_i,S_j)$ is equal to the number of arrows from the $i$-th vertex to the $j$-th vertex of $Q$ for every $S_i \in \wis{rep}_{\beta_i}~A$ and $S_j \in \wis{rep}_{\beta_j}~A$. As a consequence,
\[
Q_A = Q \qquad \text{and} \qquad \alpha_A = (1,\hdots,1) \]
In this case the associated algebra $B$ is isomorphic to $\C Q_A = \C Q$.
\end{example}

\begin{example} 
Let $A$ be a hereditary order over a smooth affine curve $X$ (or, if you prefer : a smooth Deligne-Mumford stack which is generically a curve see \cite[Coroll. 7.8]{ChanIngalls}). Then, $A$ is an affine formally smooth algebra with center $\C[X]$ and is a subalgebra of  $M_n(\C(X))$ for some $n$ (the p.i.-degree of $A$). If $\mathfrak{m}_x$ is the maximal ideal corresponding to $x$ then for all but finitely many $x$ we have that
\[
A/\mathfrak{m}_xA \simeq M_n(\C) \]
For the exceptional points $\{ x_1,\hdots,x_l \}$ (the so called ramified points) there are finitely many maximal ideals of $A$
\[
\{ P_1(i),\hdots,P_{k_i}(i) \} \quad \text{lying over} \quad \mathfrak{m}_i \]
These ideals correspond to simple representations of $A$ of dimension $\{ n_1(i),\hdots,n_{k_i}(i) \}$ such that $n_1(i) + \hdots + n_{k_i}(i) = n$.

Therefore, the representation schemes $\wis{rep}_l~A$ for $l < n$ consist of finitely many closed orbits and the semigroup generators of $\wis{comp}~A$ are given by $\alpha_j(i)$ for all $1 \leq i \leq l$ and $1 \leq j \leq k_i$ where
\[
\wis{rep}_{\alpha_j(i)}~A = \Oscr(A/P_j(i)) \rInto \wis{rep}_{n_j(i)}~A \]
It follows from \cite[Prop. 6.1]{LBlocal} that in this case the quiver $Q_A$ is the disjoint union of $l$ quivers $Q_A(i)$ of type $\tilde{A}_{k_i}$ (with circular orientation), that is,
\[
Q_A(i) = \xymatrix{
& \vtx{} \ar[r] & \vtx{} \ar[dr] & \\
\vtx{}\ar@{.}[ur] & &  & \vtx{} \ar[d]  \\
\vtx{}\ar[u] & &  & \vtx{}\ar[dl] \\
& \vtx{} \ar[ul] & \vtx{} \ar[l] & } \qquad \text{and} \qquad \alpha_A(i) = (n_1(i),\hdots,n_{k_i}(i))
\]
The associated algebra $B = B_1 \oplus \hdots \oplus B_l$ where
\[
B_i = \begin{bmatrix}
M_{n_1(i)}(\C[x]) & \hdots & M_{n_1(i) \times n_{k_i}(i)}(\C[x]) \\
\vdots & & \vdots \\
M_{n_{k_i}(i) \times n_1}(x\C[x]) & \hdots & M_{n_{k_i}(i)}(\C[x])
\end{bmatrix}
\]
(ideals below the main diagonal)
and it follows from \cite[Chpt. 9]{Reiner} or \cite[Prop. 6.1]{LBlocal} that $A$ in a neighborhood of $x_i$ is \'etale isomorphic to $B_i$.
\end{example}

\begin{example}
The modular groupalgebra $A = \C PSL_2(\Z)$ can be identified with the free algebra product
$\C \Z/2\Z \ast \C \Z/3\Z$ and therefore is a formally smooth affine algebra and every finite dimensional representation of it is isomorphic to a representation of the bipartite quiver
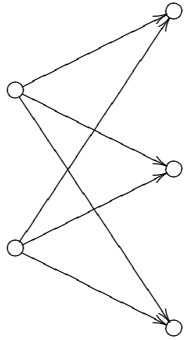
\[
\xymatrix{
& & \vtx{} \\
\vtx{} \ar[urr] \ar[drr] \ar[dddrr] & & \\
& & \vtx{} \\
\vtx{} \ar[uuurr] \ar[urr] \ar[drr]  & & \\
& & \vtx{} }
\]
where the left-hand vertex-spaces are the eigenspaces for the $\Z/2\Z$-action and those on the right-hand the eigenspaces for the $\Z/3\Z$-action. In particular, an $n$-dimensional  $\C PSL_2(\Z)$ representation has dimension vector $(a_1,a_2;b_1,b_2,b_3)$ such that $a_1+a_2 = b_1+b_2+b_3=n$ and such that the matrix determined by the $6$ arrows is invertible (it gives a base-change on the $n$-dimensional representation), see \cite{ALB} or \cite{Westbury} for more details. As a consequence $\wis{comp}~A$ has semigroup generators
\[
v_{11}=(1,0;1,0,0),~v_{12}=(1,0;0,1,0),~v_{13}=(1,0;0,0,1) \]
\[
v_{21}=(0,1;1,0,0),~v_{22}=(0,1;0,1,0),~v_{23}=(0,1;0,0,1) \]
and each component $\wis{rep}_{v_{ij}}~A$ consists of a single simple $1$-dimensional module $S_{ij}$. Because $A$ is a universal localization of the path algebra of this bipartite quiver we can compute the dimensions of $Ext^1_A(S_{ij},S_{kl})$ from the corresponding dimension vectors of the quiver. As a consequence, the associated quiver setting
$(Q_A,\alpha_A)$ is
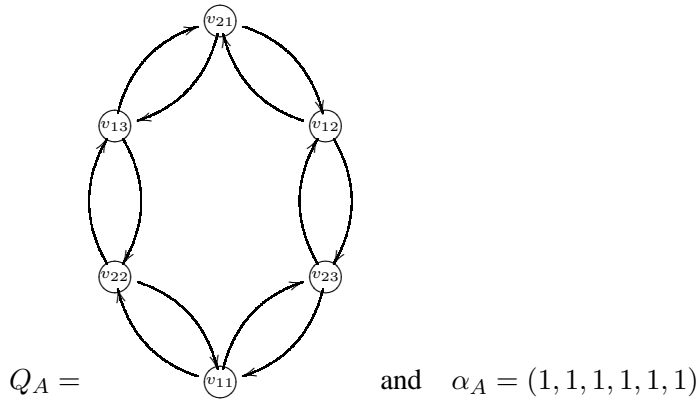
\[
Q_A = 
\xy
\POS (0,0) *\cir<6pt>{}*+{\txt\tiny{$v_{11}$}} ="S11",
\POS (14,14) *\cir<6pt>{}*+{\txt\tiny{$v_{23}$}} ="S23",
\POS (-14,14) *\cir<6pt>{}*+{\txt\tiny{$v_{22}$}} ="S22",
\POS (14,34) *\cir<6pt>{}*+{\txt\tiny{$v_{12}$}} ="S12",
\POS (-14,34) *\cir<6pt>{}*+{\txt\tiny{$v_{13}$}} ="S13",
\POS (0,48) *\cir<6pt>{}*+{\txt\tiny{$v_{21}$}} ="S21",
\POS"S11" \ar@/^2ex/ "S23"
\POS"S23" \ar@/^2ex/ "S11"
\POS"S12" \ar@/^2ex/ "S23"
\POS"S23" \ar@/^2ex/ "S12"
\POS"S12" \ar@/^2ex/ "S21"
\POS"S21" \ar@/^2ex/ "S12"
\POS"S21" \ar@/^2ex/ "S13"
\POS"S13" \ar@/^2ex/ "S21"
\POS"S13" \ar@/^2ex/ "S22"
\POS"S22" \ar@/^2ex/ "S13"
\POS"S11" \ar@/^2ex/ "S22"
\POS"S22" \ar@/^2ex/ "S11"
\endxy
 \quad \text{and} \quad \alpha_A = (1,1,1,1,1,1) \]
and the associated algebra $B$ is the path algebra $\C Q_A$.
\end{example}

\section{Simple dimension vectors}

\begin{definition} $\alpha \in \wis{comp}~A$ is said to be a simple dimension vector provided there is a non-empty Zariski open subset of $\wis{rep}_{\alpha}~A$ consisting of simple representations.
The set of all simple dimension vectors of $A$ will be denoted by $\wis{simp}~A$.
\end{definition}

\begin{theorem} If $A$ is a formally smooth algebra then $Q_A$ contains enough information to determine $\wis{simp}~A$.
\end{theorem}

\Proof
If $\alpha \in \wis{comp}~A$ we can write (possibly in several ways)
\[
\alpha = a_1 \beta_1 + \hdots + a_k \beta_k \]
with $a_i \in \N$ and $\{ \beta_1,\hdots,\beta_k \}$ the semigroup generators of $\wis{comp}~A$. This implies that there are points in $\wis{rep}_{\alpha}~A$ corresponding to semi-simple representations
\[
M = S_1^{\oplus a_1} \oplus \hdots \oplus S_k^{\oplus a_k} \]
where the $S_i$ are distinct simple representations in $\wis{rep}_{\beta_i}~A$ and we can choose the $S_i$ such that for all $1 \leq i,j \leq k$ we have that $S_i \oplus S_j \in E(\beta_i,\beta_j)$. But then the local quiver setting of $M$ in $\wis{rep}_{\alpha}~A$ is determined by $(Q_A,\epsilon)$ where $\epsilon=(a_1,\hdots,a_k)$. Because $\wis{rep}_{\alpha}~A$ is irreducible, it follows from section~1 that $\alpha \in \wis{simp}~A$ if and only if $\epsilon$ is the dimension vector of a simple representation of $Q_A$. These dimension vectors have been classified in \cite{LBProcesi} and we recall the result.

Let $\chi$ be the Euler-form of $Q_A$, that is $\chi = (c_{ij})_{i,j} \in M_k(\Z)$ with $c_{ij} = \delta_{ij} - \# \{ $ arrows from $v_i$ to $v_j~\}$ and let $\delta_i$ be the dimension vector of a vertex-simple concentrated in vertex $v_i$. Then, $\epsilon$ is the dimension vector of a simple representation of $Q_A$ if and only if the following conditions are satisfied : (1) the support $\wis{supp}(\epsilon)$ is a strongly connected subquiver of $Q_A$ (there is an oriented cycle in $\wis{supp}(\epsilon)$ containing each pair $(i,j)$ such that $\{ v_i,v_j \} \subset \wis{supp}(\epsilon)$)
and (2) for all $v_i \in \wis{supp}(\epsilon)$ we have the numerical conditions
\[
\chi(\epsilon,\delta_i) \leq 0 \qquad \text{and} \qquad \chi(\delta_i,\epsilon) \leq 0 \]
{\em unless} $\wis{supp}(\epsilon)$ in an oriented cycle of type $\tilde{A}_l$ for some $l$ in which case all components of $\epsilon$ must be equal to one.

\begin{example}
If we apply this result to the setting $(Q_A,\alpha_A)$ for the modular groupalgebra $A = \C PSL_2(\Z)$ we find that a dimension vector $(a_1,a_2;b_1,b_2,b_3)$ has an open subset of simple $PSL_2(\Z)$-representations if and only if
\[
b_j \leq a_i \qquad \text{for all $1 \leq i \leq 2$ and $1 \leq j \leq 3$} \]
which is the criterium found by Bruce Westbury in \cite{Westbury}.
\end{example}

\section{$ext(\alpha,\beta)$ on simples}

In the foregoing section we were careful to take $S_i \oplus S_j \in E(\beta_i,\beta_j)$ but this is not really necessary.

\begin{lemma} For every simple $S_i$ in $\wis{rep}_{\beta_i}~A$ and simple $S_j$ in $\wis{rep}_{\beta_j}~A$ we have that
\[
dim_{\C}~Ext^1_A(S_i,S_j) = ext(\beta_i,\beta_j) \]
\end{lemma}

\Proof
For $\alpha = \beta_i+\beta_j$ the local structure of $\wis{rep}_{\alpha}~A$ near the orbit of $M=S_i \oplus S_j$ is determined by the local quiver setting $(Q_M,\alpha_M)$ 
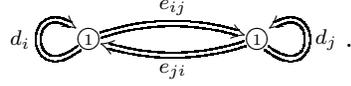
\[
\xymatrix{\vtx{1} \ar@/^/@2{->}[rr]^{e_{ij}} \ar@2@(ld, lu)^{d_i}  &  & \vtx{1}
\ar@2@(rd, ru)[]_{d_j} \ar@/^/@2{->}[ll]^{e_{ji}}}.
\]
where $d_i = dim~\wis{iss}_{\beta_i}~A$, $d_j = dim~\wis{iss}_{\beta_j}~A$, $e_{ij} = dim_{\C}~Ext^1_A(S_i,S_j)$ and $e_{ji}=dim_{\C}~Ext^1_A(S_j,S_i)$ and $\alpha_M=(1,1)$.
Now, $\wis{rep}_{\alpha_M}~Q_M$ can be identified with $Ext^1_A(M,M)$ which is the normal space $N_M$ to the orbit of $M$ in $\wis{rep}_{\alpha}~A$ which has therefore dimension
\[
dim~N_M = d_i+d_j+e_{ij}+e_{ji} \]
On the other hand, there is a point $N=S'_i \oplus S'_j$ in $\wis{rep}_{\alpha}~A$ with $N \in E(\beta_i,\beta_j)$ and $N \in E(\beta_j,\beta_i)$ and we have
\[
e_{ij} \geq ext(\beta_i,\beta_j) \qquad e_{ji} \geq ext(\beta_j,\beta_i) \]
and as for the normal space $N_N$ to the orbit of $N$ in $\wis{rep}_{\alpha}~A$ we have that
\[
dim~N_N = d_i+d_j+ext(\beta_i,\beta_j)+ext(\beta_j,\beta_i) \]
As the stabilizer subgroup of $N$ and $M$ is isomorphic to $\C^* \times \C^*$ and $\wis{rep}_{\alpha}~A$ is smooth it follows that $dim~N_M = dim~N_N$ from which the claim follows.

\begin{theorem} If $A$ is a formally smooth algebra, then $Q_A$ contains enough information to compute the dimension of $Ext^1_A(S,T)$ for all simple representations $S$ and $T$. More precisely,
\[
dim~Ext^1_A(S,T) = - \chi(\epsilon,\eta) \]
if $S$ is a simple representation in $\wis{rep}_{\alpha}~A$ where $\alpha = \sum_i a_i \beta_i$ and $\epsilon = (a_1,\hdots,a_k)$ and if $T$ is a simple representation of $\wis{rep}_{\beta}~A$ where $\beta = \sum_i b_i \beta_i$ and $\eta = (b_1,\hdots.b_k)$.
\end{theorem}

\Proof
Let $S_i$ and $T_i$ be distinct simples in $\wis{rep}_{\beta_i}~A$ and consider the semi-simple representations $M$ resp. $N$ in $\wis{rep}_{\alpha}~A$ resp. $\wis{rep}_{\beta}~A$
\[
M = S_1^{\oplus a_1} \oplus \hdots \oplus S_k^{\oplus a_k} \quad \text{and} \quad
N = T_1^{\oplus b_1} \oplus \hdots \oplus T_k^{\oplus b_k} \]
By the foregoing lemma we have complete information on the local quiver setting of $M \oplus N$
in $\wis{rep}_{\alpha+\beta}~A$ from $Q_A$. By assumption on $\alpha$ and $\beta$ there is a Zariski open subset of simples $S' \in \wis{rep}_{\alpha}~A$ and simples $T' \in \wis{rep}_{\beta}~A$ such that $S' \oplus T'$ lies in a neighborhood of $M \oplus N$. By the result of \cite{LBProcesi} we can therefore reconstruct the local quiver setting of $S' \oplus T'$ from that of $M \oplus N$. This quiver setting has the following form
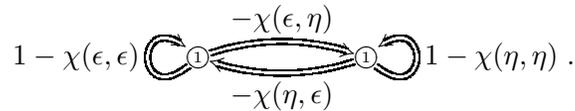
\[
\xymatrix{\vtx{1} \ar@/^/@2{->}[rr]^{\txt{$-\chi(\epsilon,\eta)$}} \ar@2@(ld, lu)^{\txt{$1-\chi(\epsilon,\epsilon)$}}  &  & \vtx{1}
\ar@2@(rd, ru)[]_{\txt{$1-\chi(\eta,\eta)$}} \ar@/^/@2{->}[ll]^{\txt{$-\chi(\eta,\epsilon)$}}}.
\]
from which we deduce that
\[
ext(\alpha.\beta) = - \chi(\epsilon,\eta) \]
But then comparing the local quiver settings of $S' \oplus T'$ with $S \oplus T$ and repeating the argument of the foregoing lemma, the result follows.

\section{The main result}

\begin{theorem} If $A$ is a formally smooth algebra, the associated quiver $Q_A$ contains enough information to reconstruct the local quiver settings $(Q_{\xi},\alpha_{\xi})$ for any semi-simple
representation
\[
M_{\xi} = S_1^{\oplus e_1} \oplus \hdots \oplus S_l^{\oplus e_l} \]
of $A$.
\end{theorem}

\Proof
This is a direct consequence of the foregoing two sections. To begin, we can determine the possible dimension vectors $\alpha_i$ of the simple components $S_i$. Write $\alpha_i = \sum_{j=1}^k a_j(i) \beta_j$ then $\epsilon_i = (a_1(i),\hdots,a_k(i))$ must be the dimension vector of a simple representation of the associated quiver $Q_A$. Moreover, by the previous theorem we know that
\[
dim~Ext^1_A(S_i,S_j) = \delta_{ij}-\chi(\epsilon_i,\epsilon_j) \]
and hence have full knowledge of the local quiver $Q_{\xi}$.

\par \vskip 4mm
Recall that the results of \cite{Bocklandt} and \cite{BLBVdW} allow us to classify the singularities of the quotient varieties $\wis{iss}_{\alpha}~A$ upto smooth equivalence and, in particular, to determine their singular loci.

\providecommand{\href}[2]{#2}
\begingroup\raggedright\endgroup


\begin{thebibliography}{10}

\bibitem{ALB}
Jan Adriaenssens and Lieven Le Bruyn, {\it Local quivers and stable representations}, Comm. Alg. (to appear)  \href{http://arxiv.org/abs/math.RT/0010251}{math.RT/0010251} (2000)

\bibitem{Bocklandt}
Raf Bocklandt, {\it Smooth quiver representation spaces}, J.Alg. {\bf 253} (2000) 296-313

\bibitem{BLBVdW}
Raf Bocklandt, Lieven Le Bruyn and Geert Van de Weyer, {\it Smooth order singularities},
J. Alg. \& Appl. (to appear) \href{http://arxiv.org/abs/math.RA/0207250}{math.RA/0207250} (2002)

\bibitem{ChanIngalls}
Daniel Chan and Colin Ingalls, {\it Noncommutative coordinate rings and stacks}, LMS (to appear)

\bibitem{CBS}
William Crawley-Boevey and Jan Schr\"oer, {\it Irreducible components of varieties of modules}, J. Reine Angew. Math. {\bf 553} (2002), 201-220

\bibitem{CuntzQuillen}
Joachim Cuntz and Daniel Quillen, {\it Algebra extensions and nonsingularity}, JAMS {\bf 8} (1995) 251-289

\bibitem{Kontsevich}
Maxim Kontsevich, {\it Formal non-commutative symplectic geometry},  Gelfand seminar 1990-1992, Birkhauser (1993) 173-187

\bibitem{LBratid}
Lieven Le Bruyn, {\it Rational identities of matrices and a theorem of G. M. Bergman}, Comm. Alg. {\bf 21(7)} (1993) 2577-2581

\bibitem{LBnagatn}
Lieven Le Bruyn, {\it noncommutative geometry@n}, \href{http://arxiv.org/abs/math.AG/9904171}{math.AG/9904171} (1999)

\bibitem{LBlocal}
Lieven Le Bruyn, {\it Local structure of Schelter-Procesi smooth orders}, Trans. AMS {\bf 352} (2000) 4815-4841

\bibitem{LBProcesi}
Lieven Le Bruyn and Claudio Procesi, {\it Semi-simple representations of quivers}, {\em Trans. AMS} {\bf 317} (1990) 585-598

\bibitem{Morrison}
Kent Morrison, {\it The connected component group of an algebra}, in {\em Lecture Notes in Mathematics} {\bf 903} (1980) 257-262

\bibitem{Procesi}
Claudio Procesi, {\it A formal inverse to the Cayley-Hamilton theorem}, J. Alg. {\bf 107} (1987) 63-74

\bibitem{Reiner}
Irving Reiner, {\it Maximal Orders}, Academic Press (1975)

\bibitem{Schofield}
Aidan Schofield, {\it Representations of rings over skew fields}, London Mathematical Society Lecture Notes Series {\bf 92} Cambridge University Press (1985)

\bibitem{Westbury}
Bruce Westbury, {\it On the character varieties of the modular group}, preprint Nottingham (1995)

\end{thebibliography}
\end{document}